\documentclass[12pt]{article}
\usepackage{graphicx}
\usepackage{color}
\usepackage{amssymb,amsmath}
\def\B{\mathcal{B}}
\def\BSP{\B_{\mathrm{sp}}}
\def\F{\mathcal{F}}
\def\FSP{\F_{\mathrm{sp}}}
\def\G{\mathcal{G}}

\def\P{\mathcal{P}}

\def\R{\mathcal{R}}
\def\S{\mathcal{S}}
\def\T{\mathcal{T}}
\def\R{\mathcal{R}}
\def\N{\mathcal{N}}
\def\M{\mathcal{M}}
\def\Aut{\mathrm{Aut}}
\newcommand{\EE}{\mathbb{E}}
\newcommand{\Leg}{\mathrm{Leg}}

\newcommand{\SPG}{\G_{\mathrm{sp}}}
\newcommand{\nn}{\mbox{\small 1} \hspace{-0,30em} 1}

\newtheorem{thm}{Theorem}
\newtheorem{corol}{Corollary}
\newtheorem{lemma}{Lemma}

\newtheorem{propos}{Proposition}

\def\proof{\noindent{\it Proof}\/.\ }

\begin{document}
\title{\bf The structure and labelled enumeration of $K_{3,3}$-subdivision-free projective-planar graphs\footnote{With the partial support of FQRNT (Qu\'ebec) and NSERC (Canada)}}
\bigbreak
\author{Andrei Gagarin, Gilbert Labelle and Pierre Leroux \\[0.1in]
 Laboratoire de Combinatoire et d'Informatique Math\'ematique,\\ 
 Universit\'e du Qu\'ebec \`a Montr\'eal, Montr\'eal, Qu\'ebec, CANADA, H3C 3P8\\[0.1in]
\small e-mail: {\texttt gagarin@lacim.uqam.ca, labelle.gilbert@uqam.ca} and \texttt leroux.pierre@uqam.ca }
\maketitle
\begin{abstract}
We consider the class $\F$ of $2$-connected non-planar $K_{3,3}$-subdivision-free graphs that are embeddable in the projective plane. We show that these graphs admit a unique decomposition as a graph $K_5$ (the \emph{core}) where the edges are replaced by two-pole networks constructed from 2-connected planar graphs. A method to enumerate these graphs in the labelled case is described. Moreover, we enumerate the homeomorphically irreducible graphs in $\F$ and homeomorphically irreducible 2-connected planar graphs. Particular use is made of two-pole (directed) series-parallel networks. We also show that the number $m$ of edges of graphs in $\F$ satisfies the bound  $m\le 3n-6$, for $n\ge 6$ vertices.
\end{abstract}
%
\section{Introduction}
\begin{figure}[h] \label{fig:K5Projective}
	\centerline {\includegraphics[height=1.6in]{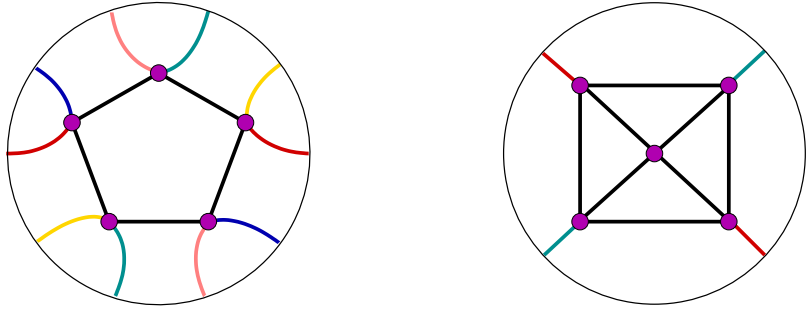}}
	\caption{\it Embeddings of $K_5$ in the projective plane.}
\end{figure}

The {\it projective plane} is a non-orientable surface of non-orientable genus $1$ that can 
be represented as a circular disk with its antipodal points identified. 
A graph $G$ is {\it projective-planar} if it can be drawn on the projective plane 
without any pair of edges crossing. 
See Figure \ref{fig:K5Projective} where two projective planar embeddings of $K_5$ are represented. 
A graph $G$ is projective-planar if and only if it contains at most one non-planar projective-planar $2$-connected component while all the other $2$-connected components of $G$ are planar.

In this paper, we consider the class $\F$ of 2-connected non-planar projective-planar graphs without a $K_{3,3}$-subdivision. 
Results for the class $C_\F$ of connected projective-planar (non-planar) graphs with no $K_{3,3}$-subdivisions are then easily deduced.
Since $K_{3,3}$ is a 3-regular graph, it is possible to see that the graphs with no $K_{3,3}$-subdivisions are precisely the graphs with no $K_{3,3}$-minors. Therefore we may refer to them as {\em $K_{3,3}$-free graphs}. By Kuratowski's Theorem \cite{Kuratowski}, these graphs must contain a subdivision of $K_5$. Hence the simplest graphs $G$ in $\F$ consist of a graph $K_5$ (the \emph{core}) where the edges are replaced by two-pole networks constructed from 2-connected planar graphs. We show that in fact this is the only possibility and moreover that the core $K_5$ of $G$ is uniquely determined as well as the two-pole networks entering in the construction. This fact is expressed in Theorem \ref{thm:FK5flecheNP} by the equation 
\begin{equation} \label{eq:K5flecheNP}
\F= K_5\uparrow \N_P.
\end{equation}
%
This property is specific to the projective plane since for other surfaces, for instance,  for the torus, more complex cores such as \emph{toroidal crowns} can occur (see \cite{GLL2,GLL3}).
 
We then apply this structure theorem to enumerate
the labelled graphs in $\F$ according to the number of vertices and edges.  
Other results include a bound on the number of edges of graphs in $\F$, which is reminescent 
of planar graphs, and the enumeration of labelled homeomorphically irreducible graphs in $\F$.

In Section $2$, we state and prove the structure theorem for the class $\F$. 
A general recursive decomposition for non-planar $K_{3,3}$-free graphs is described in Wagner \cite{Wagner} and Kelmans \cite{Kelmans}. We recall the results of Fellows and Kaschube \cite{Fellows} and  Gagarin and Kocay \cite{Us2} on the structure of non-planar graphs containing a $K_5$-subdivision of a special type 
and cite from \cite{Us2} the characterization of 2-connected non-planar $K_{3,3}$-free projective-planar graphs (the class $\F$) as graphs obtained from $K_5$ by substituting
planar networks for edges. We then prove the uniqueness of this decomposition which
establishes Theorem \ref{thm:FK5flecheNP}.

In Section 3, we first review some basic notions and terminology of labelled graphical enumeration.
The reader should have some familiarity with exponential generating functions
and their operations (addition, multiplication and composition).
See, for example, Bergeron, Labelle and Leroux \cite{BLL}, Goulden and Jackson \cite{Goulden}, 
Stanley \cite{Stanley}, or Wilf \cite{Wilf}.
We use mixed generating functions of the form 
\begin{equation}
\G(x,y)=\sum_{n\ge 0} \sum_{m\ge 0} g_{n,m}y^m\frac{x^n}{n!},
\end{equation}
where  $g_{n,m}$ is the number of graphs in a given class $\G$, with $m$ edges and on a set
of vertices $V_n$ of size $n$.
The main result here is that 
the effect on generating functions of the edge substitution operation is given by 
\begin{equation}
(\G\uparrow \N)(x,y)=\G(x,\N(x,y)).
\end{equation}
%
Use is made of the enumeration of the class $P$ of 2-connected planar graphs 
by Bender, Gao and Wormald \cite{Bender} and Bodirsky, Gr\"opl, and Kang \cite{Bodirsky}, 
based on previous work of Mullin and Schellenberg \cite{Mullin} 
on labelled enumeration of $3$-connected planar graphs.

In Section 4, we study a special class of two-pole networks, the class $\R$ of \emph{series-parallel} networks. Note that parallel edges are not permitted here. We also consider the species of series-parallel graphs denoted by $\SPG$ \cite{Duffin, Flocchini, Timothy}.  
Our presentation follows a more structural and intuitive approach, where the emphasis is put on the structure classes or species. Indeed their recursive definitions can be translated into functional equations satisfied by the species themselves and these relations are then expressed in terms of their generating functions. Moreover, many computations can be carried out and understood at the species level, before taking the generating functions.

Attention is given, in Section 5, to the enumeration of the class $H_{\F}$ 
of homeomorphically irreducible graphs in $\F$.  
Here again the edge substitution operation $H\uparrow \R$ plays a central role, 
where $\R$ represents the class of series-parallel networks. We introduce a new general iterative scheme for the computation of the generating series $H(x,y)$ satisfying an identity of the form 
\begin{equation} 
B(x,y) = H(x,R(x,y)),
\end{equation}
where $B(x,y)$ and $R(x,y)$ are known. This scheme is also applied to enumerate the class $H_{P}$ of homeomorphically irreducible 2-connected planar graphs.

Finally, a short concluding section is devoted to some related questions. For example, the class $C_\F$ of connected projective-planar (non-planar) $K_{3,3}$-free graphs is studied and asymptotic questions are touched on. Extensions to the class $\T$ of non-planar toroidal $K_{3,3}$-free graphs and to unlabelled enumeration are also considered.

Numerical results appear in six tables giving the number of labelled graphs of the 
families $\F$, $H_{\F}$ and $H_{P}$, 
for example, for $n\le16$ and $m\le42$ for the class $H_{\F}$. Results for connected graphs in $\F$ are also given.
The calculations were done with {\it Maple \hspace{-0.5mm}9.5} software on Apple Macintosh computers.

%
\section{A structure theorem for $K_{3,3}$-free projective-planar graphs} 
By convention, the graph $K_2$ is considered as a 2-connected (non-separable) graph in this paper. A \textit{two-pole network} (or more simply, a \emph{network}) is a connected graph $N$ with two distinguished vertices $0$ and $1$, such that the graph $N\cup 01$ is $2$-connected, where the notation $N\cup ab$ is used for the graph obtained from $N$ by adding the edge $ab$ if it is not already there. The vertices $0$ and $1$ are called the {\it poles} of $N$, and all the other vertices of $N$ are called \textit{internal} vertices.  

We define an operator $\tau$ acting on $2$-pole networks, $N \longmapsto \tau\cdot N$, which interchanges the poles $0$ and $1$. 
A class $\N$ of networks is called \textit{symmetric} if $N\in\N \Longrightarrow \tau\cdot N\in\N$. 
\smallskip

\noindent
{\bf Definition.}
Let $\G$ be a class of graphs and $\N$ be a symmetric class of networks.  We denote by $\G\uparrow \N$ the class of pairs of graphs $(G,G_0)$, such that
\begin{enumerate}
\item
the graph $G_0$ is in $\G$ (called the \emph{core}),
\item
the vertex set $V(G_0)$ is a subset of $V(G)$,
\item
there exists a family $\{N_e:e\in E(G_0)\}$ of networks in $\N$ (called the \emph{components}) such that the graph $G$ can be obtained from $G_0$ by substituting  $N_e$ for each edge $e\in E(G_0)$, identifying the poles of $N$ with the extremities of $e$ according to some orientation.
\end{enumerate}

\begin{figure}[h] \label{fig:exemple}
\begin{center} \includegraphics[height=2.2in]{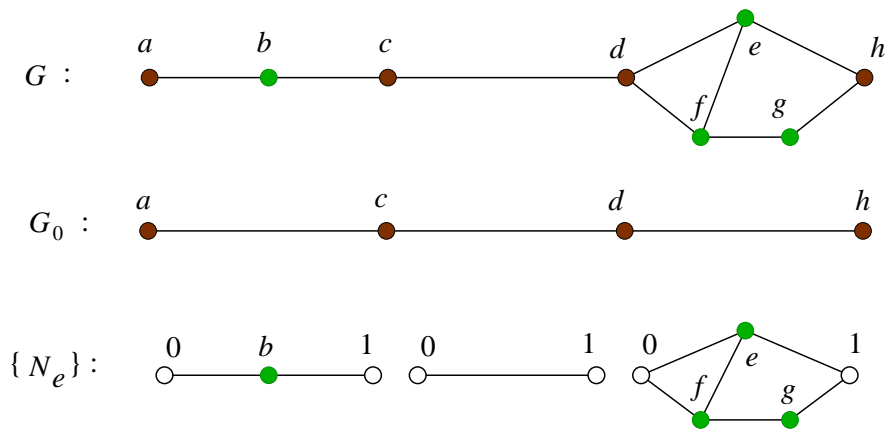}
\end{center}
\vspace{-4mm}
\caption{Example of a $(P_4 \uparrow \N)$-structure $(G,G_0)$}
\end{figure}

An example of a $(\G \uparrow \N)$-structure $(G,G_0)$, with $\G=P_4$, the class of path-graphs of order 4, and $\N$ = the class of all networks, is given in Figure \ref{fig:exemple}.

The substitution of a network $N_e$ for an edge $e$ of $G_0$ is similar to the $2$-sum operation defined for matroids and graphs by Seymour in \cite{Seymour}. One difference is that when the edge $01$ is absent from the network $N_e$, the corresponding edge $e$ is also absent from the resulting graph $G$.

We say that the composition $\G\uparrow \N$ is {\it canonical} if for any structure $(G,G_0)\in \G\uparrow \N$, the core  $G_0\in \G$ is uniquely determined by the graph $G$. In this case, we can identify $\G \uparrow \N$ with the class of resulting graphs $G$.

%

A network $N$ is {\it strongly planar} if the graph $N\cup 01$ is planar. Denote by $\N_P$ the class of strongly planar networks. Let $K_5$ denote the class of complete graphs with 5 vertices.
\begin{thm} \label{thm:FK5flecheNP}
The class $\F$ of 2-connected non-planar projective-planar $K_{3,3}$-free graphs can be expressed as a canonical composition
\begin{equation} \label{eq:FK5flecheNPbis}
\F=K_5\uparrow \N_P.
\end{equation}
\end{thm}
\proof 
We use the following previously established results. Following Diestel \cite{Diestel}, a subgraph isomorphic to a $K_5$-subdivision is denoted by $TK_5$. 
Let $G$ be a non-planar 2-connected graph with a $TK_5$.
The vertices of degree $4$ in $TK_5$ are called \emph{corners} and 
the vertices of degree $2$ are the \emph{inner vertices} of $TK_5$. 
A path connecting two corners and containing no other corner is called a 
\emph{side} of the $K_5$-subdivision. 
Note that two sides of the same $TK_5$ can have at most one common corner and no common inner vertices. 
A path $p$ in $G$ such that one endpoint is an inner vertex of $TK_5$,
the other endpoint is on a different side of $TK_5$ and all other vertices and edges 
lie in $G\backslash TK_5$ is called a \emph{shortcut} of the $K_5$-subdivision. 
A vertex $u\in G\backslash TK_5$ is called a $3$-\emph{corner vertex} with respect to $TK_5$ 
if $G\backslash TK_5$ contains internally disjoint paths from $u$ 
to at least three corners of the $K_5$-subdivision. 

The following proposition is proved in a different context: 
%
\begin{propos} [Asano \cite{Asano},\cite{Fellows,Us2}] \label{propos:3-corner} 
Let $G$ be a non-planar graph with a $K_5$-subdivision $TK_5$ for which there is either a shortcut or a $3$-corner vertex. Then $G$ contains a $K_{3,3}$-subdivision.
\end{propos}
\begin{propos} [\cite{Fellows,Us2}] \label{propos:sidecomponents}  
Let $G$ be a 2-connected graph with a $TK_5$ having neither a shortcut nor a $3$-corner vertex. 
Let $K$ denote the set of corners of $TK_5$.
Then any connected component $C$ of $G\backslash K$ contains inner vertices of at most one side 
of $TK_5$ and $C$ is adjacent to exactly two corners of $TK_5$ in $G$. 
\end{propos}
Given a graph $G$ satisfying the hypothesis of Proposition~\ref{propos:sidecomponents}, 
a \emph{side component} of $TK_5$ is defined as a subgraph of $G$ induced by a pair 
of corners $a$ and $b$ of $TK_5$ and all connected components of $G\backslash K$ which are adjacent to both $a$ and $b$. 
Notice that side components or the entire graph $G$ can contain a $K_{3,3}$-subdivision. 
For example, see Figure~3. 
However, if $G$ has no $K_{3,3}$-subdivisions, then Proposition \ref{propos:sidecomponents} 
can be applied in virtue of Proposition~\ref{propos:3-corner}.
\begin{figure}[h] \label{fig:Chrismastree}
	\centerline {\includegraphics[height=1.5in]{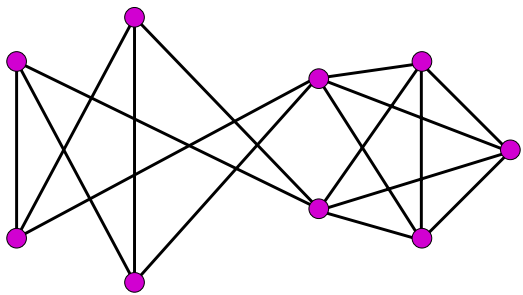}}
	\caption{\it A graph containing subdivisons of $K_{3,3}$ and $K_5$.}
\end{figure}
\begin{corol} [\cite{Fellows,Us2}] \label{corol:unique corner}
For a 2-connected graph $G$ with a $TK_5$ having no shortcut or $3$-corner vertex, 
two side components of $TK_5$ in $G$ have at most one vertex in common. 
The common vertex is the corner of intersection of two corresponding sides of $TK_5$.
\end{corol}

Thus we see that a graph $G$ satisfying the hypothesis of 
Proposition \ref{propos:sidecomponents}
can be decomposed into side components corresponding to the sides of $TK_5$.
Each side component $H$ contains exactly two corners $a$ and $b$ corresponding to a side of $TK_5$. 
If the edge $ab$ between the corners is not in $H$, we add it to $H$ to obtain $H\cup ab$. 
Otherwise $H\cup ab=H$.
We call 
$H\cup ab$ an \emph{augmented side component} of $TK_5$. 

If $G$ is a 2-connected non-planar $K_{3,3}$-free graph, we can apply 
Proposition \ref{propos:sidecomponents} and Corollary \ref{corol:unique corner} 
to decompose $G$ into side components of a $TK_5$.
Notice that with $6$ or more vertices, these graphs are not $3$-connected 
since such a $3$-connected non-planar graph contains a $K_{3,3}$-subdivision 
(see, for example, \cite{Asano}).
%
\begin{thm}[\cite{Us2}] \label{thm:projectiveplanar}
A 2-connected non-planar $K_{3,3}$-free graph $G$ is projective-planar 
if and only if $G$ contains a $K_5$-subdivision $TK_5$ 
such that each augmented side component of $TK_5$ is planar.
\end{thm}
The proof 
is based on properties of the two embeddings of the complete graph $K_5$ in the projective plane shown in Figure \ref{fig:K5Projective}.
Theorem \ref{thm:projectiveplanar} reduces projective-planarity testing to planarity testing if there is no $K_{3,3}$-subdivision in the graph. 
%
Theorem \ref{thm:projectiveplanar} can be strengthened to give the following equivalent form of Theorem \ref{thm:FK5flecheNP}. We say that a side component $H$ is \emph{strongly planar} if the augmented side component $H\cup ab$ is planar.  This is coherent with the previously defined concept of strongly planar network. 
\begin{thm}  \label{thm:unicity}
A 2-connected non-planar projective-planar $K_{3,3}$-free graph $G$ has a unique decomposition into  strongly planar side components of a $TK_5$.
\end{thm}
\proof 
By Proposition \ref{propos:sidecomponents}, the set of corners $K$ of $TK_5$ completely 
defines a decomposition into the side components. 
Therefore it is sufficient to show that any other $K_5$-subdivision $TK_5^\prime$ in $G$ shares the same set $K$ of corners with $TK_5$.
Since each augmented side component of $TK_5$ in $G$ is planar, all corners of $TK_5^\prime$ cannot be contained in any particular side component. 
Suppose that a corner $a$ of $TK_5^\prime$ is not in $K$. 
Then  $a$ is in a side component $S$ of $TK_5$. 
Recall that there should be four disjoint paths from $a$ to the four other corners of $TK_5^\prime$.
Since there is no shortcut or $3$-corner vertex of $TK_5$ in $G$, 
the side component $S$ of $TK_5$ must contain at least $2$ other corners of $TK_5^\prime$, say $b$ and $c$. Now consider a corner $d$ of $TK_5^\prime$ that is not in the side component $S$ of $TK_5$. 
Three of the sides adjacent to $d$ must connect $d$ to the three corners 
$a, b$ and $c$ of $TK_5^\prime$. 
However the disjoint sides $da, db$ and $dc$ of $TK_5^\prime$ must share 
two corners of the side component $S$ of $TK_5$, a contradiction. \phantom{i}
\hfill\rule{2mm}{2mm}

This concludes the proof of Theorem \ref{thm:FK5flecheNP}.
\hfill\rule{2mm}{2mm}

\medskip
A corollary to Euler's formula for the plane says that a planar graph with $n\ge 3$ vertices 
can have at most $3n-6$ edges (see, for example, 
\cite{Diestel}). Let us state this for 2-connected planar graphs with $n$ vertices and $m$ edges as follows:
\begin{equation} \label{eq:mplanar}
m\le  \left\{ \begin{array}{ll} 
 3n-5 & \mbox{if \,$n=2$} \\ 
 3n-6 & \mbox{if \,$n\ge 3$} 
\end{array}.
\right. 
\end{equation} 
In fact, $m=3n-5=1$ if $n=2$. The generalized Euler Formula (see, for example, \cite{Carsten}) implies 
that a projective-planar graph $G$ with $n$ vertices can have up to $3n-3$ edges. 
An arbitrary $K_{3,3}$-free graph $G$ is known to have at most $3n-5$ edges 
(see \cite{Asano}). 
However we show here that projective-planar $K_{3,3}$-free graphs satisfy the following stronger relation, 
which is similar to that of planar graphs. 

\begin{propos}  \label{propos:mprojplanar}
The number $m$ of edges of a non-planar projective-planar $K_{3,3}$-free $n$-vertex graph $G$ satisfies $m=3n-5=10$ if $n=5$, and
\begin{equation} \label{eq:mprojplanar}
m\le 3n-6 \ \mbox{if}\  n\ge 6.
\end{equation}
\end{propos}
\proof It is sufficient to prove the result for  2-connected graphs. By Theorem $2$, each augmented side component $S_i$ of $G$, $i=1,2,\ldots,10$, 
satisfies the condition $(\ref{eq:mplanar})$
with $n=n_i$, the number of vertices and $m=m_i$, the number of edges of $S_i$, $i=1,2,\ldots,10$. 
Since each corner of $TK_5$ is in precisely $4$ side components, we have 
$\sum_{i=1}^{10}n_i=n+ 15$ and we obtain, by summing these 10 inequalities, 
$$
m=\sum_{i=1}^{10}m_i\le \left\{ \begin{array}{lll} 
3\sum_{i=1}^{10}n_i-50 = 3(n+15)-50 = 3n-5 & \mbox{if $n=5$} \\ \\
3\sum_{i=1}^{10}n_i-51 = 3(n+15)-51 = 3n-6 & \mbox{if $n\ge 6$} 
\end{array},
\right.
$$
since 
$n=5$ iff $n_i=2$, $i=1,2,\ldots,10$, and $n\ge 6$ if and only if at least one $n_j\ge 3$, $j=1,2,\ldots,10$.
\hfill\rule{2mm}{2mm}

\medskip
Notice that Corollary 8.3.5 of \cite{Diestel} implies that graphs 
without a $K_5$-subdivision also can have at most $3n-6$ edges.

%
\section{Initial enumerative results}
We now consider the labelled enumeration of projective-planar $K_{3,3}$-free graphs according to the numbers of vertices and edges.
We first review some basic notions and terminology of labelled enumeration.
The reader should have some familiarity with exponential generating functions
and their operations (addition, multiplication and composition).
See \cite{BLL}, \cite{Goulden}, \cite{Stanley}, or \cite{Wilf}.

By a \emph{labelled} graph, we mean a simple graph $G=(V,E)$ where the set of vertices $V=V(G)$ 
is itself the set of labels and the labelling function is the identity function.
$V$ is called the \emph{underlying set} of $G$.
An edge $e$ of $G$ then consists of an unordered pair $e=uv$ of elements of $V$ 
and $E=E(G)$ denotes the set of edges of $G$.
If $W$ is another set and $\sigma:V\tilde{\rightarrow}W$ is a bijection, 
then any graph $G=(V,E)$ on $V$ 
can be transformed into a graph $G^\prime=\sigma(G)=(W,\sigma(E))$, where
$\sigma(E)=\{\sigma(e) =\sigma(u)\sigma(v)\,|\, e\in E\}$. 
We say that $G^\prime$ is obtained from $G$ by \emph{vertex relabelling}
and that $\sigma$ is a graph \emph{isomorphism} $G\tilde{\rightarrow}G^\prime$.
An \emph{unlabelled graph} is then seen as an isomophism class $\gamma$ of labelled graphs. 
We write $\gamma=\gamma(G)$ if $\gamma$ is the isomorphism class of $G$. 
By the \emph{number of ways to label} an unlabelled graph $\gamma(G)$, 
where $G=(V,E)$, we mean the number of distinct graphs $G^\prime$ on 
the underlying set $V$ which are isomorphic to $G$. 
Recall that this number is given by $n!/|\Aut(G)|$, 
where $n=|V|$ and $\Aut(G)$ denotes the automorphism
group of $G$.

A \emph{species} of graphs is a class of labelled graphs which is closed under vertex relabellings. 
Thus any class $\G$ of unlabelled graphs gives rise to a species, also denoted by $\G$,
by taking the set union of the isomorphism classes in $\G$.
For any species $\G$ of graphs, we introduce its \emph{mixed (exponential) generating function}
$\G(x,y)$ as the formal power series
\begin{equation} \label{eq:Gdexy}
\G(x,y)=\sum_{n\ge 0}g_n(y)\frac{x^n}{n!}, \ \ 
\mbox {with}\ \ g_n(y)=\sum_{m\geq0} g_{n,m}y^m,
\end{equation}
where $g_{n,m}$ is the number of graphs in $\G$ with $m$ edges over a given set
of vertices $V_n$ of size $n$. Here $y$ is a formal variable which acts as an edge counter.
For example, for the species $\G=K=\{K_n\}_{n\geq0}$ of complete graphs, we have
\begin{equation} \label{eq:Kxy}
K(x,y)=\sum_{n\geq0}y^{n\choose2}x^n/n!,
\end{equation}
while for the species $\G=\G_a$ of all simple graphs, we have $\G_a(x,y)=K(x,1+y)$. 
Another example is the class of \emph{discrete} graphs (i.e. with no edges),
which we denote by $\EE$ (for French "Ensemble") since these are just sets of vertices, 
and we have $\EE(x,y)=\sum_{n\geq0}x^n/n!=\exp(x)$.  

A species of graphs is \emph{molecular} if it contains only one isomorphism class. Examples include the class $K_1$ of one-vertex graphs, which is denoted by $X$ and satisfies $X(x,y)=x$, and the class $K_2$, with $K_2(x,y)=yx^2/2$.
In general, for a molecular species $\gamma=\gamma(G)$, where $G$ has $n$ vertices and $m$ edges, we have $\gamma(x,y)= \frac{y^mn!}{|\Aut(G)|}x^n/n!=y^mx^n/|\Aut(G)|$. For example, we have 
\begin{equation} \label{eq:K5xy}
K_5(x,y)= x^5y^{10}/5!
\end{equation}

\medskip
For two-pole networks, only the internal vertices 
form the underlying set for the purpose of enumeration and for species considerations. In particular, the mixed generating function of a class (or species) $\N$ of networks is defined by 
\begin{equation} \label{eq:Ndexy}
\N(x,y)=\sum_{n\ge 0}\nu_n(y)\frac{x^n}{n!}, \ \ 
\mbox {with}\ \ \nu_n(y)=\sum_{m\geq0} \nu_{n,m}y^m,
\end{equation}
where $\nu_{n,m}$ is the number of networks in $\N$ with $m$ edges with a given set of internal vertices $V_n$ of size $n$. 

\medskip
There is an operator $\tau$ acting on two-pole networks, $N \longmapsto \tau\cdot N$, which interchanges the poles $0$ and $1$. A species $\N$ of networks is called \textit{symmetric} if $N\in\N \Longrightarrow \tau\cdot N\in\N$. Examples of symmetric species of networks include the class $\N_P$ of strongly planar networks and the class $\R$ of series-parallel networks described in the next section.

\begin{propos} [T. Walsh \cite{Timothy}] \label{lemme:GflecheN}
Let $\G$ be a species of graphs and $\N$ be a symmetric species of networks.
Then we have 
\begin{equation} \label{eq:GflecheNxy}
(\G\uparrow \N)(x,y)=\G(x,\N(x,y)).
\end{equation}
\end{propos}
\proof  
If $(G,G_0)$ is a $(\G\uparrow \N)$-structure where the core graph $G_0$ has $k$ edges, thus contributing a term $y^k$ to $\G(x,y)$, we can assume that the underlying set of $G_0$ is linearly ordered. We say that the substitution of a network $N$ for an edge $e=ab$, with $a<b$, is \emph{coherent} if the pole 0 of $N$ is identified with $a$ and the pole 1, with $b$.
Since the class $\N$ is symmetric, it is sufficient to restrict ourselves to coherent substitutions. Moreover, we can order the edges of $G_0$ lexicographically so that the process of edge subtitution is uniquely determined by a list of $k$ disjoint networks in $\N$. 
Since these lists are counted by $\N^k(x,y)$, formula (\ref{eq:GflecheNxy}) follows. 
\hfill\rule{2mm}{2mm}
\begin{corol} 
The mixed generating function $\F(x,y)$  of labelled 2-connected non-planar projective-planar $K_{3,3}$-free graphs is given by
\begin{equation} \label{eq:Fxy}
\F(x,y)=\frac{x^5\N_P^{10}(x,y)}{5!}.
\end{equation}
\end{corol}
\proof This follows from Theorem \ref{thm:FK5flecheNP}, Proposition \ref{lemme:GflecheN}, 
and the fact that $K_5(x,y)= x^5y^{10}/5!$. 
\hfill\rule{2mm}{2mm}

\medskip
There remains to compute the generating series $\N_P(x,y)$.
If $\M$ is a class of networks which do not contain the edge $01$, then we denote by $y\,\M$ the class obtained by adding this edge to all the networks
of $\M$. Observe that there are two distinct networks on the empty set, 
namely the trivial network $\nn$, consisting of two isolated poles 0 and 1, and the one edge network $y\nn$. 

Now let $B$ be a given species of $2$-connected graphs 
containing $K_2$, 
for example $B=B_a$, the class of all $2$-connected graphs, $B=\{K_2\}$
or, more importantly here, $B=P$, the class of all $2$-connected \emph{planar} graphs.
We denote by $B^{(y)}$ the species of graphs obtained by selecting and removing an edge
in all possible ways from graphs in $B$ . 
Note that 
\begin{equation}
B^{(y)}(x,y)=\frac{\partial}{\partial y}B(x,y).
\end{equation}
If, moreover, the endpoints of the selected edge are unlabelled and numbered 0 and 1, 
in all possible ways, the resulting class of networks is denoted by $B_{0,1}$.
Relabelling the two poles yields the identity 
\begin{equation}
x^2B_{0,1}(x,y)=2\,B^{(y)}(x,y).
\end{equation}
%
Finally, we introduce the species of networks $\N_B$ associated to the class $B$ 
by the formula
\begin{equation} \label{eq:NB}
\N_B = B_{0,1} + yB_{0,1} - \nn=(1+y)B_{0,1} - \nn.
\end{equation}
Thus, the generating function of $\N_B$ is given by
\begin{equation} \label{eq:NBxy}
\N_B(x,y) = (1+y)\frac{2}{x^2}\frac{\partial}{\partial y}B(x,y) - 1.
\end{equation}
%

%
\begin{table}[!h] \label{table:fnm}
\centerline {\scriptsize
	\begin{tabular}[t]{|| r | r | r || r | r | r || r | r | r ||}
	\hline
	$n$ & $m$ & $f_{n,m}$ & $n$ & $m$ & $f_{n,m}$ & $n$ & $m$ & $f_{n,m}$\\
	\hline \hline
	5 & 10 & 1 &  11 & 16 & 1664863200 & 13 & 18 & 1261490630400\\
	\hline
	6 & 11 & 60 & 11 & 17 & 17556739200 & 13 & 19 & 21330659750400\\
	\hline
	6 & 12 & 60 & 11 & 18 & 78956539200 & 13 & 20 & 159781461840000\\
	\hline
	7 & 12 & 2310 & 11 & 19 & 202084621200 & 13 & 21 & 713882464495200\\
	\hline
	7 & 13 & 5250 & 11 & 20 & 334016949420 & 13 & 22 & 2168012582255520\\
	\hline
	7 & 14 & 3150 & 11 & 21 & 387489624060 & 13 & 23 & 4841896937557680\\
	\hline
	7 & 15 & 210 & 11 & 22 & 335202677040 & 13 & 24 & 8367745313108610\\
	\hline
	8 & 13 & 73920 & 11 & 23 & 221055080400 & 13 & 25 & 11501380415300490\\
	\hline
	8 & 14 & 283920 & 11 & 24 & 107529691500 & 13 & 26 & 12648862825333020\\
	\hline
	8 & 15 & 380240 & 11 & 25 & 35726852700 & 13 & 27 & 11024998506341820\\
	\hline
	8 & 16 & 205520 & 11 & 26 & 7205814000 & 13 & 28 & 7476620617155690\\
	\hline
	8 & 17 & 40320 & 11 & 27 & 663616800 & 13 & 29 & 3846042558007650\\
	\hline
	8 & 18 & 5040 & 12 & 17 & 45664819200 & 13 & 30 & 1446703666808400\\
	\hline
	9 & 14 & 2162160 & 12 & 18 & 617512896000 & 13 & 31 & 374735495534400\\
	\hline
	9 & 15 & 12383280 & 12 & 19 & 3642195110400 & 13 & 32 & 59680805184000\\
	\hline
	9 & 16 & 27592740 & 12 & 20 & 12576897194400 & 13 & 33 & 4401725328000\\
	\hline
	9 & 17 & 30616740 & 12 & 21 & 28943910959040 & 14 & 19 & 35321737651200\\
	\hline
	9 & 18 & 18419940 & 12 & 22 & 48122268218640 & 14 & 20 & 732123289497600\\
	\hline
	9 & 19 & 6656580 & 12 & 23 & 61023477279600 & 14 & 21 & 6797952466905600\\
	\hline
	9 & 20 & 1678320 & 12 & 24 & 60601323301200 & 14 & 22 & 38137563765100800\\
	\hline
	9 & 21 & 196560 & 12 & 25 & 46937904829200 & 14 & 23 & 147357768378300480\\
	\hline
	10 & 15 & 60540480 & 12 & 26 & 27584940398400 & 14 & 24 & 423597368531216880\\
	\hline
	10 & 16 & 481572000 & 12 & 27 & 11793019392000 & 14 & 25 & 951297908961680280\\
	\hline
	10 & 17 & 1578301200 & 12 & 28 & 3448102996800 & 14 & 26 & 1715806516686001740\\
	\hline
	10 & 18 & 2810039400 & 12 & 29 & 615367368000 & 14 & 27 & 2511869870973763300 \\
	\hline
	10 & 19 & 3055603320 & 12 & 30 & 50494752000 & 14 & 28 & 2981167142609535880\\
	\hline
	10 & 20 & 2214739800 & \multicolumn{3}{c ||}{ } & 14 & 29 & 2845977828319866240\\
	\cline {1-3} \cline {7-9}
	10 & 21 & 1155735000 &  \multicolumn{3}{c ||}{ } & 14 & 30 & 2159624129854611420\\
	\cline {1-3} \cline {7-9}
	10 & 22 & 432356400 &  \multicolumn{3}{c ||}{ } & 14 & 31 & 1281613625914642020\\
	\cline {1-3} \cline {7-9}
	10 & 23 & 98809200 &  \multicolumn{3}{c ||}{ } & 14 & 32 & 580974136160418000\\
	\cline {1-3} \cline {7-9}
	10 & 24 & 10281600 &  \multicolumn{3}{c ||}{ } & 14 & 33 & 194019911828542800\\
	\cline {1-3} \cline {7-9}
	\multicolumn{6}{c ||}{ } & 14 & 34 & 44947147269024000\\
	\cline {7-9}
	\multicolumn{6}{c ||}{ } & 14 & 35 & 6446992644892800\\
	\cline {7-9}
	\multicolumn{6}{c ||}{ } & 14 & 36 & 431053060838400\\
	\cline {7-9}
	\end{tabular} }
	\caption{The number $f_{n,m}$ of labelled non-planar projective-planar 2-connected graphs without 	a $K_{3,3}$-sub\-di\-vi\-sion (having $n$ vertices and $m$ edges).}
\end{table}

\smallskip
Let $P$ denote the species of $2$-connected planar graphs. Then the associated class $\N_P$ of networks described above is precisely the class of strongly planar networks.
%
Methods for computing the generating function $P(x,y)$ of labelled $2$-connected planar graphs are described in \cite{Bender} and \cite{Bodirsky}.
Both methods are based on the network decomposition of \cite{Trakh} 
which is also stated for planar graph embeddings in \cite{Tutte}.  
The decomposition allows to count the $2$-connected planar graphs via labelled $3$-connected planar graphs whose counting can be derived from \cite{Mullin}.

Formulas (\ref{eq:NBxy}) and (\ref{eq:Fxy}) can then be used to compute $\N_P(x,y)$ and $\F(x,y)$. Numerical results are presented in Tables 1 and 2,
where $\F(x,y)=\sum_{n\geq5}\sum_{m} f_{n,m}x^ny^m/n!$ and $f_n=\sum_{m}f_{n,m}$.
%
\begin{table}[h] \label{table:fn}
\centerline {\scriptsize
	\begin{tabular}[t]{|| r | r ||}
	\hline
	$n$ & $f_n$ \\
	\hline \hline
	5 & 1 \\
	\hline
	6 & 120\\
	\hline
	7 & 10920\\
	\hline
	8 & 988960\\
	\hline
	9 & 99706320\\
	\hline
	10 & 11897978400\\
	\hline
	11 & 1729153068720\\
	\hline
	12 & 306003079514880\\
	\hline
	13 & 64657337524631280\\
	\hline
	14 & 15890834362452489440\\
	\hline
	15 & 4435396700216405763840\\
	\hline
	16 & 1379778057502074926142720\\
	\hline
	17 & 471689356958791639787042560\\
	\hline
	18 & 175335742043846629500183667200\\
	\hline
	19 & 70291642269058321415718042668160\\
	\hline
	20 & 30195035473057938652243866755197440\\
	\hline 
	\end{tabular} }
	\caption{The number $f_n$ of labelled non-planar projective-planar 2-connected graphs without a $K_{3,3}$-sub\-di\-vi\-sion (having $n$ vertices).}
\end{table}
%
\section{Series-parallel networks and graphs}
In this section, we study a special class of two-pole networks, the class $\R$ of \emph{series-parallel} networks (also called two-terminal directed series-parallel networks). Note that parallel edges are not permitted here. We also consider the species of series-parallel graphs, denoted by $\SPG$. See, for example, \cite{Duffin}, \cite{Flocchini} and \cite{Timothy}.  

It is assumed that the poles of a network $N$ are distinct from those of any other network. 
There are two main operations on two-pole networks: parallel composition and series composition. Let $S$ be a finite set of disjoint networks which are not equal to $\nn$ and do not contain the edge $01$. The \emph{parallel composition} of $S$ is the network obtained by taking the union of the graphs in $S$ where, moreover, all the $0$-poles are fusioned into one $0$-pole and similarly for the $1$-poles. By convention, the parallel composition of an empty set of networks is the trivial network $\nn$.
If $\N$ is a species of networks which are distinct from $\nn$ and have non-adjacent poles 
and if each network in a class $\M$ can be viewed unambiguously as a parallel
composition of networks in $\N$, then the result can be expressed as a species composition $\M=\EE(\N)$, and we have 
\begin{equation} \label{eq:parallelxy}
\M(x,y)=\exp(\N(x,y)).
\end{equation}
Note that the class $\N$ is then included in $\M$ and that $\nn$ is in $\M$. 

Let $M$ and $N$ be two non-trivial disjoint networks. 
The \emph{series composition} $M\cdot_sN$ of $M$ followed by $N$ is a network
whose underlying set is the union of the underlying sets of $M$ and $N$ plus an extra element.
It is obtained by 
taking the graph union of $M$ and $N$ where moreover the $1$-pole of $M$ is fusioned
with the $0$-pole of $N$ and this \emph{connecting vertex} is labelled by the extra element.
The \emph{series composition} $\M\cdot_s\N$ of two species of networks $\M$ and $\N$ 
not containing $\nn$ 
is the class obtained by taking all series compositions $M\cdot_sN$ 
with $M\in\M$ and $N\in\N$.
If moreover the two components $M\in\M$ and $N\in\N$ are uniquely determined by the resulting
network $M\cdot_sN$, the species $\M\cdot_s\N$ can be expressed as the species product 
$\M X\N$, where the factor $X$ corresponds to the connecting vertex, and we have 
\begin{equation} \label{eq:series}
(\M\cdot_s\N)(x,y)=x\M(x,y)\N(x,y).
\end{equation}
%

The species $\R$ of \emph{series-parallel} networks can be defined recursively as  
the smallest class of networks containing the one-edge network $y\nn$ and closed under
series and parallel compositions.
We partition $\R$ as $\R=\S+\P$, where $\S$ represents the species of \emph{essentially series} networks (i.e. of the form (\ref{eq:seriesparallel}) below) 
and $\P=\R-\S$ is the complementary class, of \emph{essentially parallel} networks. 
These classes are characterized recursively by the following functional equations 
involving series and parallel composition: 
\begin{equation} \label{eq:seriesparallel}
\S = \P\cdot_s\R = \P X\R.
\end{equation}
\begin{equation} \label{eq:parallelseries}
\R=(1+y)\EE(\S)-\nn.
\end{equation}
From these two equations, we deduce 
\begin{equation} \label{eq:PdeR}
\R = \S+\P = \P X \R+\P = \P(1+X\R) \ \Rightarrow \ \P = \R/(1+X\R)
\end{equation}
and 
\begin{equation} \label{eq:RdeR}
\R = (1+y)\EE(\P X\R)-\nn = (1+y)\EE(\frac{X\R^2}{1+X\R})-\nn
\end{equation}
and for the generating functions,
\begin{equation} \label{eq:RxdeRx}
\R(x,y) = (1+y)\exp\left(\frac{x\R^2(x,y)}{1+x\R(x,y)}\right)-1.
\end{equation}
The series $\R(x,y)$ can be computed recursively using (\ref{eq:RxdeRx}).

\medskip
Now  a \emph{series-parallel graph} is a 2-connected graph which is either an edge or can be obtained from a series-parallel network by adding the edge $01$ (if not already present) and labelling the poles. Series-parallel graphs can be characterized as 2-connected graphs without a $K_4$-subdivision (see \cite{Duffin}). The class of series-parallel graphs is denoted by $\SPG$. It is easy to see that the networks induced by series-parallel graphs are precisely the series-parallel networks, i.e. that
\begin{equation} \label{eq:NSPG}
\N_{\SPG} = \R.
\end{equation}
This implies, using (\ref{eq:NBxy}), that
\begin{equation} \label{eq:SPGxy}
\SPG(x,y) = \frac{x^2}{2}\int_{0}^{y} \frac{\R(x,t)+1}{1+t} \,dt. 
\end{equation}
%

\section{Homeomorphically irreducible graphs}
A graph is called \emph{homeomorphically irreducible} if it contains no vertex of degree 2. For a graph $G$ embeddable in a surface, any subdivision of $G$ is trivially embeddable in the same surface. Therefore, it is interesting to count graphs embeddable in a surface that are minimal with respect to the operation of subdivision, i.e. homeomorphically irreducible graphs. Here we do this for the classes $P$ of $2$-connected planar graphs and $\F$ of 2-connected non-planar projective-planar $K_{3,3}$-free graphs, applying the method of Walsh (\cite{Timothy}) as follows.

Any 2-connected graph $G$ is either a series-parallel graph or contains a unique 
2-connected homeomorphically irreducible core $C(G)$, which is different from $K_2$, and unique components $\{N_e\}_{e\in E(C(G))}$ which are series-parallel networks, whose composition gives $G$.
Let $\B$ be a species of 2-connected graphs. Denote by $H_{\B}$ the class of graphs which are homeomorphically irreducible cores of graphs in $\B$. Also set $\BSP = \B\cap\SPG$ which is the class of series-parallel graphs in $\B$. 
\begin{propos} \label{prop:BBSHBR}
Let $\B$ be a species of 2-connected graphs such that 
\begin{enumerate}
\item
$H_{\B}$ is contained in $\B$,
\item
$\B$ is closed under edge substitution by series-parallel networks, 
i.e. $\B\uparrow \R$ is contained in $\B$.
\end{enumerate}
Then we have 
\begin{equation} \label{eq:BBSHBR}
\B = \BSP + H_{\B}\uparrow \R,
\end{equation}
and the composition $H_{\B}\uparrow \R$ is canonical.
\end{propos}

For the generating functions, it follows that
\begin{equation} \label{eq:BBSHBRxy}
\B(x,y) = \BSP(x,y) + H_{\B}(x,\R(x,y)),
\end{equation}
from which the series $H_{\B}(x,y)$ can be computed, in virtue of the following lemma:
\begin{lemma} \label{lemma:inversion}
Let $B(x,y)$ and $R(x,y)$ be two-variable formal power series such that $R(x,y)=y+O(y^2)$.
Then there exists a unique formal power series $H(x,y)$ such that 
\begin{equation} \label{eq:BHxRxy}
B(x,y) = H(x,R(x,y)).
\end{equation}
Moreover, $H(x,y)$ can be expressed as 
\begin{equation} \label{eq:Hxydelta}
H(x,y) = \sum_{i\geq0}(-1)^i\Delta_R^iB(x,y),
\end{equation}
where $\Delta_R$ is an operator defined on two-variable formal power series $F(x,y)$ by 
$$\Delta_RF(x,y)=F(x,R(x,y))-F(x,y).$$
\end{lemma}
\proof The first statement follows from the fact that under the hypothesis, the series  
$R(x,y)$, viewed as a formal power series in the variable $y$, is invertible under composition.
The equation (\ref{eq:Hxydelta}) then follows from the observation that equation (\ref{eq:BHxRxy})
is equivalent to $B(x,y)=(I+\Delta_R)H(x,y)$, where $I$ denotes the identity operator. 
Details are left to the reader. 
\hfill\rule{2mm}{2mm}
\begin{table}[!h]
\centerline {\scriptsize
	\begin{tabular}[t]{|| r | r | r || r | r | r || r | r | r ||}
	\hline
	$n$ & $m$ & $H_P(n,m)$ & $n$ & $m$ & $H_P(n,m)$ & $n$ & $m$ & $H_P(n,m)$\\
	\hline \hline
	2 & 1 & 1 &  10 & 15 & 5700240 & 13 & 20 & 1845922478400\\
	\hline
	4 & 6 & 1 & 10 & 16 & 297561600 & 13 & 21 & 74599125400800\\
	\hline
	5 & 8 & 15 & 10 & 17 & 2930596200 & 13 & 22 & 989130437895600\\
	\hline
	5 & 9 & 10 & 10 & 18 & 12343659300 & 13 & 23 & 6630351423696000\\
	\hline
	6 & 9 & 60 & 10 & 19 & 28301918400 & 13 & 24 & 26817549328369800\\
	\hline
	6 & 10 & 477 & 10 & 20 & 38982967065 & 13 & 25 & 71682957811565100\\
	\hline
	6 & 11 & 585 & 10 & 21 & 33331061925 & 13 & 26 & 133187371098982200\\
	\hline
	6 & 12 & 195 & 10 & 22 & 17392158000 & 13 & 27 & 176696593868094300\\
	\hline
	7 & 11 & 4410 & 10 & 23 & 5088258000 & 13 & 28 & 169059691482031350\\
	\hline
	7 & 12 & 23520 & 10 & 24 & 641277000 & 13 & 29 & 116043129855402750\\
	\hline
	7 & 13 & 37800 & 11 & 17 & 1659042000 & 13 & 30 & 55840786515914400\\
	\hline
	7 & 14 & 24570 & 11 & 18 & 41399542800 & 13 & 31 & 17912090131135200\\
	\hline
	7 & 15 & 5712 & 11 & 19 & 340745605200 & 13 & 32 & 3443842956153600\\
	\hline
	8 & 12 & 13440 & 11 & 20 & 1407085287300 & 13 & 33 & 300495408595200\\
	\hline
	8 & 13 & 332640 & 11 & 21 & 3435723903150 & 14 & 21 & 4217639025600\\
	\hline
	8 & 14 & 1543860 & 11 & 22 & 5355953687700 & 14 & 22 & 594554129769600\\
	\hline
	8 & 15 & 2917740 & 11 & 23 & 5504170275450 & 14 & 23 & 15443454480854400\\
	\hline
	8 & 16 & 2708160 & 11 & 24 & 3728003340600 & 14 & 24 & 172534400373535200\\
	\hline
	8 & 17 & 1237320 & 11 & 25 & 1606084131960 & 14 & 25 & 1084459459555672200\\
	\hline
	8 & 18 & 223440 & 11 & 26 & 399801679200 & 14 & 26 & 4361354314691635800\\
	\hline
	9 & 14 & 2177280 & 11 & 27 & 43859692800 & 14 & 27 & 12056896085921586900\\
	\hline
	9 & 15 & 28962360 & 12 & 18 & 3996669600 & 14 & 28 & 23900258922899609250\\
	\hline
	9 & 16 & 126168840 & 12 & 19 & 362978431200 & 14 & 29 & 34803857129521580100\\
	\hline
	9 & 17 & 266535360 & 12 & 20 & 6150939628680 & 14 & 30 & 37649908340175226095\\
	\hline
	9 & 18 & 311551380 & 12 & 21 & 44916513919200 & 14 & 31 & 30254152933093434345\\
	\hline
	9 & 19 & 207170460 & 12 & 22 & 183180611357100 & 14 & 32 & 17843305708519691400\\
	\hline
	9 & 20 & 73710000 & 12 & 23 & 470167225050600 & 14 & 33 & 7512030324951352200\\
	\hline 
	9 & 21 & 10929600 &  12 & 24 & 807689258734050 & 14 & 34 & 2139154225643635200\\
	\hline 
	\multicolumn{3}{c ||}{ } &  12 & 25 & 956591057815470 & 14 & 35 & 369529809669820800\\
	\cline {4-9}
	\multicolumn{3}{c ||}{ } &  12 & 26 & 786477564207720 & 14 & 36 & 29262949937020800\\
	\cline {4-9}
	\multicolumn{3}{c ||}{ } &  12 & 27 & 442142453075400 & \multicolumn{3}{c}{ }\\
	\cline {4-6}
	\multicolumn{3}{c ||}{ } & 12 & 28 & 162493688649600 & \multicolumn{3}{c}{ }\\
	\cline {4-6}
	\multicolumn{3}{c ||}{ } & 12 & 29 & 35240506963200 & \multicolumn{3}{c}{ }\\
	\cline {4-6}
	\multicolumn{3}{c ||}{ } & 12 & 30 & 3424685806080 & \multicolumn{3}{c}{ }\\
	\cline {4-6}
	\end{tabular} }
	\caption{ The number $H_P(n,m)$ of labelled 2-connected homeomorphically irreducible planar graphs (having $n$ vertices and $m$ edges).}
\end{table}
%

%
\begin{table}[t]
\centerline {\scriptsize
	\begin{tabular}[t]{|| r | r ||}
	\hline
	$n$ & $H_P(n)$ \\
	\hline \hline
	4 & 1\\
	\hline
	5 & 25\\
	\hline
	6 & 1317\\
	\hline
	7 & 96012\\
	\hline
	8 & 8976600\\
	\hline
	9 & 1027205280\\
	\hline
	10 & 139315157730\\
	\hline
	11 & 21864486188160\\
	\hline
	12 & 3898841480307900\\
	\hline
	13 & 778680435365714700\\
	\hline
	14 & 172192746831203449890\\
	\hline
	15 & 41765231538761743574100\\
	\hline
	16 & 11024455369912310561835600\\
	\hline
	17 & 3146065407516184280981053200\\
	\hline 
	18 & 965135197612755256313598822450\\
	\hline 
	19 & 316731891055609655106993297185400\\
	\hline 
	20 & 110718818921232836033343337842628500\\
	\hline 
	\end{tabular} }
	\caption{ The number $H_P(n)$ of labelled 2-connected homeomorphically irreducible planar graphs (having $n$ vertices).}
\end{table}
%
\begin{table}[!h]
\centerline {\scriptsize
	\begin{tabular}[t]{|| r | r | r || r | r | r || r | r | r ||}
	\hline
	$n$ & $m$ & $h_{n,m}$ & $n$ & $m$ & $h_{n,m}$ & $n$ & $m$ & $h_{n,m}$\\
	\hline \hline
	5 & 10 & 1 &  12 & 21 & 2025777600 & 15 & 25 & 7205830632000\\
	\hline
	7 & 14 & 210 & 12 & 22 & 44347564800 & 15 & 26 & 923081887728000\\
	\hline
	7 & 15 & 210 & 12 & 23 & 321609657600 & 15 & 27 & 21992072494392000\\
	\hline
	8 & 15 & 3360 & 12 & 24 & 1163155593600 & 15 & 28 & 226159100164998000\\
	\hline
	8 & 16 & 13440 & 12 & 25 & 2450459088000 & 15 & 29 & 1307760868616202000\\
	\hline
	8 & 17 & 15120 & 12 & 26 & 3214825059600 & 15 & 30 & 4819747766658224400\\
	\hline
	8 & 18 & 5040 & 12 & 27 & 2674115413200 & 15 & 31 & 12125014783013632500\\
	\hline
	9 & 16 & 15120 & 12 & 28 & 1375742491200 & 15 & 32 & 21647164674205570500\\
	\hline
	9 & 17 & 257040 & 12 & 29 & 400355524800 & 15 & 33 & 27986100453182371500\\
	\hline
	9 & 18 & 948780 & 12 & 30 & 50494752000 & 15 & 34 & 26352994794020744850\\
	\hline
	9 & 19 & 1372140 & 13 & 22 & 5772967200 & 15 & 35 & 17930879200055571750\\
	\hline
	9 & 20 & 861840 & 13 & 23 & 462940077600 & 15 & 36 & 8598757529558196000\\
	\hline
	9 & 21 & 196560 & 13 & 24 & 7019020008000 & 15 & 37 & 2759781730918212000\\
	\hline
	10 & 18 & 2116800 & 13 & 25 & 45946008108000 & 15 & 38 & 532536328868544000\\
	\hline
	10 & 19 & 23511600 & 13 & 26 & 167127038278200 & 15 & 39 & 46746961252992000\\
	\hline
	10 & 20 & 85428000 & 13 & 27 & 378396335116800 & 16 & 27 & 5038469339904000\\
	\hline
	10 & 21 & 145681200 & 13 & 28 & 563392705525650 & 16 & 28 & 277876008393984000\\
	\hline
	10 & 22 & 128898000 & 13 & 29 & 563307338043450 & 16 & 29 & 5018906168115840000\\
	\hline
	10 & 23 & 57531600 & 13 & 30 & 375765580990800 & 16 & 30 & 45917601694892928000\\
	\hline
	10 & 24 & 10281600 & 13 & 31 & 160783623000000 & 16 & 31 & 254980573605117360000\\
	\hline
	11 & 19 & 6652800 & 13 & 32 & 39987851904000 & 16 & 32 & 945355953679641504000\\
	\hline
	11 & 20 & 301039200 & 13 & 33 & 4401725328000 & 16 & 33 & 2474074000472282208000\\
	\hline
	11 & 21 & 2559249000 & 14 & 24 & 2746116172800 & 16 & 34 & 4725628918340842512000\\
	\hline
	11 & 22 & 9235749600 & 14 & 25 & 100222343020800 & 16 & 35 & 6713043612043772203200\\
	\hline
	11 & 23 & 17753412600 & 14 & 26 & 1207921516401600 & 16 & 36 & 7147503225193821890400\\
	\hline
	11 & 24 & 19736016300 & 14 & 27 & 7362246043152000 & 16 & 37 & 5690265588000873079200\\
	\hline
	11 & 25 & 12781345500 & 14 & 28 & 26958084641888400 & 16 & 38 & 3341472585422235264000\\
	\hline
	11 & 26 & 4491471600 &  14 & 29 & 64675702745854200 & 16 & 39 & 1406015363067771456000\\
	\hline
	11 & 27 & 663616800 &  14 & 30 & 106440372932493600 & 16 & 40 & 401357157281144064000\\
	\hline
	\multicolumn{3}{c ||}{ } &  14 & 31 & 122731243349715000 & 16 & 41 & 69663390944539392000\\
	\cline {4-9}
	\multicolumn{3}{c ||}{ }&  14 & 32 & 99327282369915600  & 16 & 42 & 5553155150839296000\\
	\cline {4-9}
	\multicolumn{3}{c ||}{ }&  14 & 33 & 55396720246467600  & \multicolumn{3}{c}{ }\\
	\cline {4-6}
	\multicolumn{3}{c ||}{ }&  14 & 34 & 20312069633856000  & \multicolumn{3}{c}{ }\\
	\cline {4-6}
	\multicolumn{3}{c ||}{ }&  14 & 35 & 4413395543356800  & \multicolumn{3}{c}{ }\\
	\cline {4-6}
	\multicolumn{3}{c ||}{ }&  14 & 36 & 431053060838400  & \multicolumn{3}{c}{ }\\
	\cline {4-6}
	\end{tabular} }
	\caption{The number $h_{n,m}$ of labelled non-planar projective-planar 2-connected $K_{3,3}$-free homeomorphically irreducible graphs 
	(having $n$ vertices and $m$ edges).}
\end{table}

Walsh uses Proposition \ref{prop:BBSHBR}  
to enumerate all labelled homeomorphically irreducible 2-connected graphs in \cite{Timothy}.
Here we give two other applications. 
First, we take $\B=P$, the class of $2$-connected planar graphs. 
In this case, $H_{\B} = H_P$ is the class of $2$-connected planar graphs 
with no vertices of degree less than $3$ and $H_P + K_2$ is the class of 2-connected
homeomorphically irreducible planar graphs.
Series-parallel graphs are known to be planar: they do not contain a subdivision of $K_4$, 
but $K_5$ and $K_{3,3}$ do.
It follows that $P_S = P\cap\SPG = \SPG$.
It is clear that the hypotheses of Proposition \ref{prop:BBSHBR} are satisfied and
we deduce from (\ref{eq:BBSHBR}) that 
\begin{equation} \label{eq:PSPG}
P = \SPG + H_P\uparrow \R
\end{equation}
where the composition is canonical. Taking the generating functions, we obtain: 
%
\begin{propos} \label{prop:PHPxy}
The mixed generating functions of the species $P$ of planar 2-con\-nected graphs
and $H_P$ of homeomorphically irreducible graphs in $P$ are related by the equation
\begin{equation} \label{eq:PHPxy}
P(x,y) =  \SPG(x,y) + H_P(x,\R(x,y)). 
\end{equation}
%
\end{propos}

  We have used this equation to compute the first terms of the series 
$H_P(x,y)=\sum_{n\geq4}\sum_{m} H_P(n,m)x^ny^m/n!$ using Lemma \ref{lemma:inversion}. The results are presented in Tables 3 and 4,
where $H_P(n) = \sum_{m}H_P(n,m)$.
Notice that the computational results of Table~3 in comparison to those of \cite{Bender} 
verify that any maximal planar graph with $n\ge 4$ vertices (and $3n-6$ edges) 
have all vertex degrees at least $3$. In other words, for $n\ge 4$ and $m=3n-6$,
we have $H_P(n,m)=P(n,m)$.

  Now we take $\B = \F$, the species of non-planar projective-planar 2-connected graphs without 
a $K_{3,3}$-subdivision. Then $H_\F$ is the class of homeomorphically irreducible graphs
in $\F$ and $\FSP=\F\cap\SPG$ is empty. 
It is clear that the hypotheses of Proposition \ref{prop:BBSHBR} are satisfied and
we deduce from (\ref{eq:BBSHBR}) that 
\begin{equation} \label{eq:FSPG}
\F = H_\F\uparrow \R
\end{equation}
where the composition is canonical. Taking the generating functions, we obtain:
\begin{propos} \label{prop:FHFxy}
The generating functions $\F(x,y)$ and $H_\F(x,y)$ of labelled non-planar projective-planar
2-connected $K_{3,3}$-free graphs and those with no vertices 
of degree less than $3$ (resp.) are related by the equation
\begin{equation} \label{eq:FHFxy}
\F(x,y) =  H_\F(x,\R(x,y)).
\end{equation}
\end{propos}

  We have used this equation to compute the first terms of the series 
$H_\F(x,y)=\sum_{n\geq5}\sum_{m} h_{n,m}x^ny^m/n!$. The results are presented in Tables 5 and 6, where $h_n = \sum_{m}h_{n,m}$.
 Notice that numbers in Table~6 are much smaller than corresponding numbers in Table~2. 
However, for $n\ge 7$ and $m=3n-6$, the corresponding numbers in Tables~1 and~5 are the same.  
This  verifies that maximal non-planar projective-planar $K_{3,3}$-free graphs have vertex degrees at least $3$. That can be seen as a corollary  to Theorem $2$, Proposition $4$ and the corresponding statement for maximal planar graphs.

%
\begin{table}[h]
\centerline {\scriptsize
	\begin{tabular}[t]{|| r | r ||}
	\hline
	$n$ & $h_n$ \\
	\hline \hline
	5 & 1 \\
	\hline
	6 & 0\\
	\hline
	7 & 420\\
	\hline
	8 & 36960\\
	\hline
	9 & 3651480\\
	\hline
	10 & 453448800\\
	\hline
	11 & 67528553400\\
	\hline
	12 & 11697130922400\\
	\hline
	13 & 2306595939347700\\
	\hline
	14 & 509359060543132800\\
	\hline
	15 & 124356566550728011500\\
	\hline
	16 & 33226132945543622884800\\
	\hline
	17 & 9635384706205021006042800\\
	\hline
	18 & 3012126613564117021370798400\\
	\hline
	19 & 1009263543337906919842715535600\\
	\hline
	20 & 360698621436519186180018210552000\\
	\hline 
	\end{tabular} }
	\caption{The number of labelled non-planar projective-planar 2-connected $K_{3,3}$-free homeomorphically irreducible graphs 
	(having $n$ vertices).}
\end{table}

\medskip
There is an alternate way to compute the series $H_\F(x,y)$ which reduces computations significantly. The idea is to determine what are the side components that should be substituted into the edges of $K_5$ in order to obtain homeomorphically irreducible graphs in $\F$. The question is to determine the class of networks $\N$ for which $H_\F=K_5\uparrow \N$. 
A first try 
is to take the class $\N=\N_{H_P}$ of networks $N$ such that $N\cup01$
is a planar 2-connected graph with no vertices of degree less than 3. Following 
(\ref{eq:NBxy}), we have
\begin{equation} \label{eq:NHPxy}
\N_{H_P}(x,y)=(1+y)\frac{2}{x^2}\frac{\partial}{\partial y}H_P(x,y)-1,
\end{equation}
where $H_P(x,y)$ is given by Proposition \ref{prop:PHPxy}.
However, the degree requirements can be relaxed for the two poles of the network $N$ since they will become identified with the corners of $K_5$ and will have degree at least $4$.
In particular, the one-edge network $y\nn$ can be used as a side component.
Another way to have a pole of degree one is to start with a network $N$ of $\N=\N_{H_P}$ 
and \emph{add a leg} at the 0-pole, the 1-pole, or both. This means taking the series compositions
$y\nn\cdot_sN$, $N\cdot_sy\nn$, or $y\nn\cdot_sN\cdot_sy\nn$. 
In this way, we obtain the class 
$y\nn X\N_{H_P} + \N_{H_P}Xy\nn +  y\nn X\N_{H_P}Xy\nn = (2y\nn X + (y\nn X)^2)\N_{H_P}$.
Let us denote by $\Leg$ this operator of adding legs, with $\Leg(x,y)=2yx+y^2x^2$.
Now it is possible to join the 
poles by an edge to obtain one or two poles of degree 2, 
giving rise to the operator $y\Leg$, and to iterate this process. 
Hence we set
\begin{equation} \label{eq:NLegNHP}
\N_{\ell} = (1+\Leg)\left(\sum_{k\geq0}(y\Leg)^k\right)\N_{H_P}.
\end{equation}
The generating series of this class of networks is given by 
\begin{equation} \label{eq:NLegNHPxy}
\N_{\ell}(x,y)=\frac{1+\Leg(x,y)}{1-y\Leg(x,y)}\N_{H_P}(x,y),
\end{equation}
where $\N_{H_P}(x,y)$ is defined by (\ref{eq:NHPxy}), and we have the following proposition.
\begin{propos} \label{prop:HFK5flecheN}
Let $\N_{\ell}$ be the species of networks defined by equation \emph{(\ref{eq:NLegNHP})}. 
Then the species $H_\F$ of homeomorphically irreducible non-planar projective-planar 2-connected $K_{3,3}$-free graphs can be expressed as 
\emph{
\begin{equation} \label{eq:HFK5flecheN}
H_\F=K_5\uparrow(\N_{\ell} + y\nn)
\end{equation}
}
and its generating series satisfies
\begin{equation} \label{eq:HFK5flecheNxy}
H_\F(x,y)= \frac{x^5(\N_{\ell}(x,y)+y)^{10}}{5!},
\end{equation}
where $\N_{\ell}(x,y)$ is given by \emph{(\ref{eq:NLegNHPxy})}.
\end{propos}
%

%
\section{Concluding Remarks}
We have obtained the labelled enumeration of the class $\F$ of 2-connected non-planar projective-planar $K_{3,3}$-free graphs. In this short section, we mention some extensions which can be obtained of these results.

\subsection{Connected graphs in $\F$}

It is easy to deduce the labelled enumeration for the class $C_\F$ of 1-connected (i.e. connected) non-planar projective-planar $K_{3,3}$-free graphs. Indeed it suffices to attach arbitrary vertex-rooted connected planar graphs at each vertex of graphs in $\F$ in order to obtain all graphs in $C_\F$. More precisely, we have 
\begin{equation} \label{eq:CPpointxy}
C_\F(x,y)=\F(C_P^{\bullet}(x,y),y),
\end{equation}
where $C_P^{\bullet}$ denotes the class of vertex rooted connected planar graphs.

Recall that $P$ denotes the class of 2-connected planar graphs. 
Then it is well known (see \cite{BLL}, \cite{Noy}) that 
\begin{equation} \label{eq:CPpointxy}
C_P^{\bullet}(x,y)=x\exp(P^\prime(C_P^{\bullet}(x,y),y)),
\end{equation}
where $P^\prime(x,y)=\frac{\partial}{\partial x}P(x,y)$, from which $C_P^{\bullet}(x,y)$ and then $C_\F(x,y)$ can be computed. For example, setting $y=1$, we obtain the following first numbers $|C_\F[n]|$ of labelled connected non-planar projective-planar $K_{3,3}$-free graphs with $n$ vertices, $n=5,6,\dots, 20$:

\medskip
\noindent
{\small 
1, 150, 16800, 1809360, 206725050, 26484163020,  
3942600552660, 694822388340960, \\ 
145100505844928205,  35439528882292735200, \\
9927470411345581984890, 3128005716477250367216640, \\
1090689073286188397027568380, 415560636438834909293721364320, \\
 171338083303545263513720985887520, 75873636257232699557453120820157440.
}

\medskip

\subsection{Asymptotics}
Using results of Bender, Gao and Wormald \cite{Bender} and Gim\'enez and Noy \cite{Noy}, it is easy to see that labelled non-planar $K_{3,3}$-free projective-planar 2-connected and 1-connected graphs share similar asymptotic behaviours, in particular, the same growth constants, as their planar counterparts.

\subsection{Other extensions}

Using methods of the present paper, we have been able to give a similar characterization for the class $\T$ of non-planar $K_{3,3}$-free toroidal 2-connected graphs. In this case, more complex \emph{toroidal} cores can occur, whose class is denoted by $\T_C$, for which the equation 
$\T=\T_C\uparrow\N_P$ holds. See \cite{GLL2} for details.

Walsh in \cite{Timothyunlabelled} has shown how to enumerate \emph{unlabelled} graphs in a class which admits an unambiguous representation of the form $\G\uparrow \N$.
Therefore the characterization of Theorem~\ref{thm:FK5flecheNP} also leads to the  unlabelled enumeration of $K_{3,3}$-free projective-planar and toroidal graphs. This has been carried out in the paper~\cite{GLL3}.

%
\medbreak
\subsection*{ Acknowledgement.} 
The authors are thankful to Professor Timothy Walsh 
for providing useful references and discussions.

\end{document}